\documentclass[11pt]{article}
\usepackage{color}
\usepackage{amsmath,amsthm,amssymb,comment}
\usepackage{amsfonts} %
\newtheorem{theo}{Theorem}[section]
\newtheorem{lem}[theo]{Lemma}

\newcommand{\mysection}[1]{\section{#1} \setcounter{equation}{0}}
\newcommand{\dint}{\displaystyle \int}

\newcommand{\R}{\mathbb{R}}

\newcommand{\be}{\begin{equation} \label}
\newcommand{\ee}{\end{equation}}
\newcommand{\bes}{\begin{equation} \begin{array}{c} \label}
\newcommand{\ees}{\end{array} \end{equation}}
\newcommand{\bea}{\begin{eqnarray}\label}
\newcommand{\eea}{\end{eqnarray}}
\newcommand{\beas}{\begin{eqnarray} \begin{array}{rcl} \label}
\newcommand{\eeas}{\end{array} \end{eqnarray}}
\newcommand{\bas}{\begin{eqnarray*}}\newcommand{\eas}{\end{eqnarray*}}
\newcommand{\bass}{\begin{eqnarray*} \begin{array}{rcl}}
\newcommand{\eass}{\end{array} \end{eqnarray*}}
\newcommand{\basss}{\begin{eqnarray*} \begin{array}{c}}
\newcommand{\easss}{\end{array} \end{eqnarray*}}
\newcommand{\bit}{\begin{itemize}}
\newcommand{\eit}{\end{itemize}}
\newcommand{\nn}{\nonumber}
\newcommand{\eps}{\varepsilon}
\newcommand{\abs}{\\[3mm]}
\newcommand{\parab}{{\cal{P}}}

\newcommand{\AD}{{\cal{A}}_D}
\newcommand{\uv}{\underline{v}}
\newcommand{\ov}{\overline{v}}
\newcommand{\ren}{\mathbb{R}^n}
\begin{document}
\title{\bf Rate of Convergence to Barenblatt Profiles \\ for the Fast
Diffusion Equation}
\author{\Large
Marek Fila\footnote{fila@fmph.uniba.sk}  \\[4pt]
\Large Juan Luis V\'azquez\footnote{juanluis.vazquez@uam.es} \\[4pt]
\Large Michael Winkler\footnote{michael.winkler@math.uni-paderborn.de} \\[4pt]
\Large Eiji Yanagida\footnote{yanagida@math.titech.ac.jp}
}
\date{}
\maketitle

\begin{abstract} \noindent We study the asymptotic behaviour of positive
solutions of the Cauchy problem
for the fast diffusion equation as $t$ approaches the extinction time. We
find a continuum of rates of convergence to a self-similar profile.
These rates depend explicitly on the spatial decay rates of initial data.

\

\end{abstract}

\mysection{Introduction}

We consider the Cauchy problem for the fast diffusion equation
\begin{equation}\label{fd}
      \left\{ \begin{array}{ll}
      u_\tau = \nabla\cdot \left(u^{m-1}\nabla u\right),
              \qquad & y\in\R^n, \ \tau\in (0,T), \\[2mm]
      u(y,0)=u_0(y)\ge 0, \qquad & y\in\R^n,
      \end{array} \right.
\end{equation}
where $m<1$, $n>2$ and $T>0$.  It is known that for $m<m_c:=(n-2)/n$ all
solutions with initial data in some convenient space, like $L^p(\ren)$ with
$p=n(1-m)/2$, extinguish in finite time. We will always work in this range,
$m<m_c$, and consider solutions which vanish in a finite time $\tau=T$.
The purpose of this paper is to study the behaviour of such solutions near
extinction.

The situation is simpler and illustrative for our purposes in the exponent range $m>m_c$
that includes both fast ($m<1$) and slow ($m>1$) diffusion. Then solutions do not extinguish in finite time and the description of the asymptotic behaviour of  the global-in-time
solutions of (\ref{fd}) as $\tau\to\infty$ for $m> m_c$ is a
very active subject. For nonnegative data $u_0\in L^1(\R^n)$ it is proved
(cf. \cite{Vazbeh}) that all
solutions converge up to scaling to one of the so-called Barenblatt solutions, precisely
to the one with the same total mass, after proper rescaling of the solutions. For $m_c<m<1$
the Barenblatt solutions take the self-similar form
\[
U_{D}(y,\tau):=(\tau+a)^{-\beta n} \left(D+\frac{\beta(1-m)}{2}
\left|\frac{y}{(\tau+a)^\beta}\right|^2\right)^{-\frac{1}{1-m}},\qquad
\beta:= \frac{1}{n\,|m_c-m|}\, .
\]
In other words, the  Barenblatt
solutions with exponent $m>m_c$ play the role of the Gaussian solution of the linear diffusion equation in describing the asymptotic behaviour of a very wide class of nonnegative
solutions, i.e., those with integrable initial data, cf. \cite{Vbook}.  The entropy
method has allowed to provide rates for that convergence, cf.
\cite{CJMTU, CLMT, CV, DM, DT}.

This study has been extended recently to the behaviour near
extinction for $m\leq m_c$. Indeed, for $m<m_c$ and the problem posed in the whole space,
the book \cite{Vsmooth} contains a general description of the
phenomenon of extinction. It is explained there that the
occurrence of extinction depends on the size of initial data,
and  also that  different initial data may give rise to  different
extinction rates, even for the same extinction time; this may happen
for all $0<m<m_c$.  It is also proved in \cite{Vsmooth} and references
quoted there that the size of the initial data at infinity (the tail
of $u_0$) is very important in determining both the extinction time
and the decay rates.

Now, some natural questions arise, like: Are there solutions that continue the Barenblatt
formulae for the extinction range $m<m_c$; if yes, how attractive  are they?  In fact, 
those solutions exist and have the self-similar form
\begin{equation}\label{baren.form1}
U_{D, T}(y,\tau):=\frac 1{R(\tau)^n} \left(D+\frac{\beta(1-m)}{2}
\left|\frac{y}{R(\tau)}\right|^2\right)^{-\frac{1}{1-m}},
\end{equation}
where for $m<m_c$ we put $R(\tau):=(T-\tau)^{-\beta},$ and
\[
\beta= \frac1{n(1-m)-2}=\frac{1}{n\,(m_c-m)}>0.
\]

Many papers (\cite{BBDGV, BDGV, BGV, DM}, for example)
are concerned with the stabilization as $\tau\to T$ of solutions
of (\ref{fd}). The general attractive character of the Barenblatt solutions of the range $m>m_c$ is lost but
still they have a {\sl basin of attraction} formed by solutions with data that are close to the Barenblatt initial data in some norm. The study of such question is taken up in the papers
\cite{BBDGV, BDGV}, which establish convergence with  rates  of the solutions of
(\ref{fd}) towards a unique attracting limit state in that family.
More precisely, the decay rates of
\[
R(\tau)^{n}(u(\tau,y)-U_{D, T}(y,\tau))
\]
as $\tau\to T$ are discussed there (note that $R(\tau)^{n}U$ has size 1).

The  critical exponent
\[
m_*:=\frac{n-4}{n-2}<m_c\, ,
\]
has the property that the difference of two generalized Barenblatt
solutions is integrable for $m\in (m_*,m_c)$, while it is not integrable for
$m\leq m_*$. The exponent $m_*$ plays a very important role in the
results of \cite{BBDGV, BDGV, BGV}. The proofs of convergence with rates are based on
the study of the decay in time of a certain {\sl relative entropy} and a careful analysis
of the linearized problem which leads to certain functional inequalities of 
 Hardy-Poincar\'e type. In particular, the basin of attraction of a function $U_{D,T}$ in the range $m<m_*$ contains functions $u_0$ such that \[
U_{D_1,T}(\cdot,0)\leq u_0\leq U_{D_2,T}(\cdot,0),\qquad
|u_0(x)-U_{D,T}(x,0)|\in L^1(R^n).
\]
We call this set the variational basin, and for this the entropy method gives precise decay rates (the variational rates).

However we had the feeling that the basin of attraction is larger if we allow ourselves to get out of the variational framework,  where the differences of solutions do not have finite relative entropy anymore. The analysis of this possibility was done in a very interesting limit case that occurs if we take $D=0$ in formula
(\ref{baren.form1}), and we find the singular solution
\[
U_{0, T}(y,\tau):=k_*\,(T-\tau)^{\mu/2}|y|^{-\mu},
\qquad k_*:=(2(n-\mu))^{\mu/2},\qquad \mu:=\frac{2}{1-m}\,.
\]
Its attractivity properties have been studied in \cite{FVW}, confirming our guess in the form
of a continuum of slow convergence rates for increasingly larger deviations at infinity with 
 respect to the tail of $U_{0, T}$, i.\,e., $C |y|^{-\mu}$.


The purpose of the present paper is to perform the construction of solutions
in the basin of attraction of the generalized Barenblatt solutions $U_{D, T}$ with both $D,T>0$ in the exponent range $m<m_*$. Our solutions fall out of the already mentioned variational basin and consequently have a slower rate of convergence.

To study the asymptotic profile as $\tau\to T$,
it is convenient to rescale the flow and rewrite
(\ref{fd}) in self-similar variables by introducing  the time-dependent
change of variables
\begin{equation}\label{eq:chgvariable}
t:=\frac{1}{\mu}\log\left(\frac{R(\tau)}{R(0)}\right)\quad\mbox{and}\quad
x:=\sqrt{\frac{\beta}{\mu}}\,\frac y{R(\tau)}\,,
\end{equation}
with $R$ as above, and the rescaled function
\begin{equation}\label{eq:chgvariable2}
v(x,t):= R(\tau)^{n}\,u(y,\tau).
\end{equation}
In these new variables, the generalized Barenblatt functions $U_{D,T}
(y,\tau)$ are transformed into \emph{generalized Barenblatt profiles} $V_D(x)$, which
have the advantage of being stationary:
\begin{equation}\label{newBaren}
V_D(x):=(D+|x|^2)^{-1/(1-m)},\quad x\in \R^n\,.
\end{equation}
The convergence theorems in the new variables take the form of stabilization to non-trivial equilibria.  If $u$ is a solution to~(\ref{fd}), then $v$
solves the {\sl rescaled fast diffusion equation}
\begin{equation}\label{FPeqn}
v_t=\nabla\cdot \left(v^{m-1}\nabla v\right)+\mu \,\nabla\cdot(x\,v), \quad t> 0\,,\quad x\in \R^n\, ,
\end{equation}
which is a nonlinear Fokker-Planck equation. We put
as initial condition $v_0(x):=R(0)^{-n}\,u_0(y)$, where $x$ and~$y$ are
related according to (\ref{eq:chgvariable}) with $\tau=0$.
We have taken the
precise form of this transformation from \cite{BBDGV}. Note also that the
factor $\mu$ in equation (\ref{FPeqn}) can be eliminated by manipulating
the change
of variables, but then the expression of the Barenblatt solutions would
contain new
constants.

By $n >2$ and $m< m_*$, we have  $\mu+2<n$ so that the number
\[
l_\star:=\frac{n+\mu+2}{2}
\]
satisfies $\mu+2<l_\star<n$. Note that
as $m\to m_*$ we have $\mu\to n-2$ and $l_*\to n$. 

\input epsf
\vskip6mm
\setbox1=\hbox{\epsfysize 5.5truecm\epsfbox[25 235 570 560]{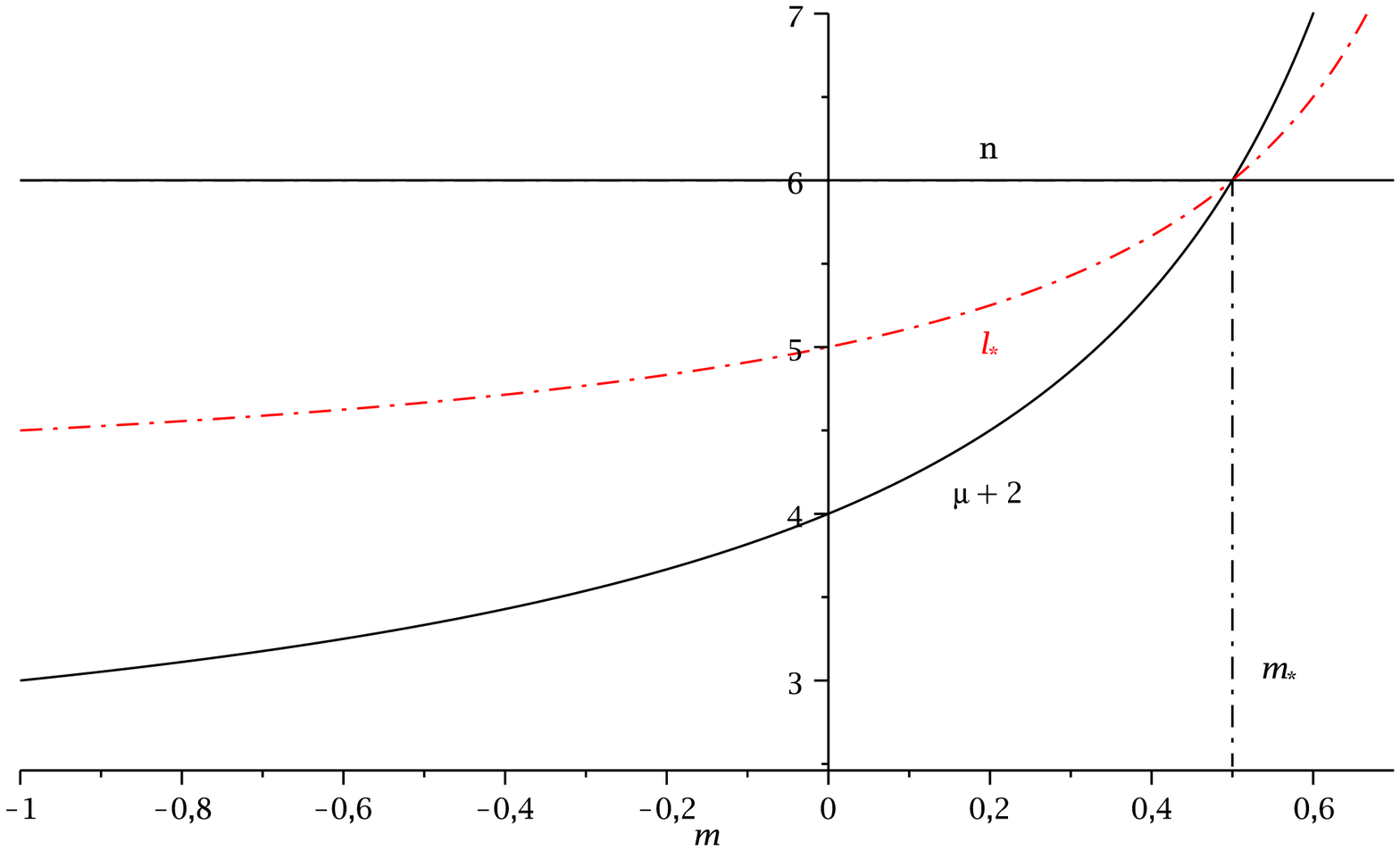}}
\vtop{
\centerline{\copy1}
\vskip3mm
\centerline{Figure 1. Values of $l_*$ and $\mu+2$ as functions of $m$ for $n=6$}
}
\vskip3mm plus3mm

Our main result is the following:
\begin{theo}\label{mainthm}
Let $m<m_*$, $n>2$. Assume that $c, D >0$ and
$
\mu + 2 < l < l_\star$.
\item[\rm(i)] If
\begin{equation}
      |v_0(x) - V_D(x)| \leq c |x|^{-l}, \qquad |x|\geq 1,
\end{equation}
and if
$$ 0<v_0(x)\leq V_\delta(x), \qquad x\in\R^n
$$
for some $\delta <D$,  then there exists $C_1>0$ such that the solution $v$ of {\rm (\ref{FPeqn})}
with the initial condition
\begin{equation} \label{eq:in}
v(x,0)=v_0(x), \qquad x \in \R^n,
\end{equation}
 satisfies
\begin{equation}
     \sup_{x\in\R^n} |v(x,t) - V_D(x)| \leq C_1 e^{-(l-\mu-2)(n-l)\,t},
      \qquad t\ge 0\,,
\end{equation}
where $(l-\mu-2)(n-l)>0$.

\item[\rm(ii)] If
\begin{equation}
      v_0(x) \leq V_D(x) - c |x|^{-l}, \qquad |x|\geq 1,
\end{equation}
and
\[
0<v_0(x)\leq V_D(x),\qquad x\in\R^n,
\]
then there exists $C_2>0$ such that the solution $v$ of {\rm
(\ref{FPeqn})},
\eqref{eq:in}
satisfies
\begin{equation}
     \sup_{x\in\R^n} \Big(V_D(x)-v(x,t)\Big) \geq C_2 e^{-(l-\mu-2)(n-l)\,t},
      \qquad t\ge 0\,.
\end{equation}

\item[\rm(iii)] If
\begin{equation}
      v_0(x) \geq V_D(x) + c |x|^{-l}, \qquad |x|\geq 1,
\end{equation}
and
\[
v_0(x)\geq V_D(x),\qquad x\in\R^n,
\]
then there exists $C_3>0$ such that the solution $v$ of {\rm
(\ref{FPeqn})}, \eqref{eq:in}
satisfies
\begin{equation}
     \sup_{x\in\R^n} \Big(v(x,t) - V_D(x)\Big) \geq C_3
e^{-(l-\mu-2)(n-l)\, t},
      \qquad t\ge 0.
\end{equation}
\end{theo}

\noindent {\bf Remarks.} (1) First of all, the result gives a
sharp description of the basin of attraction of generalized Barenblatt
solutions for $m<m_*$. For such $m$ it was shown
in \cite{BBDGV, BDGV, DS} that the basin of attraction
of a generalized Barenblatt solution $U_{D,T}$ contains all solutions
corresponding to data which, besides being trapped between two Barenblatt
profiles $U_{D_0,T}$, $U_{D_1,T}$ for the same value of $T$,
are integrable perturbations of $U_{D,T}$. Theorem~\ref{mainthm}~(i) implies
that
if $0<u_0(y)\leq U_{\delta,T}(y,0)$, $\delta >0$, and
\[
|u_0(y)-U_{D,T}(y,0)|\leq c|y|^{-l},\qquad |y|\geq 1
\]
for some $c>0$ and $l>\mu +2$ then the corresponding solution is contained
in the basin of attraction of $U_{D,T}$. Since $m<m_*$, it follows that
$\mu +2<n$. Hence, non-integrable perturbations of $U_{D,T}$ may still yield
convergence to $U_{D,T}$. The condition $l>\mu +2$ is optimal since the
difference of two Barenblatt solutions is of the order $|y|^{-(\mu +2)}$.

\noindent (2) We have thus found a continuum of convergence rates which depend explicitly on the tail of initial data. The rate $(l-\mu-2)(n-l)$ converges to zero as $l\to \mu +2$ and to 
the  maximum value
\be{alpha_star}
	\alpha_\star:=\frac{(n-\mu-2)^2}{4}
\ee
as $l\to l_\star$. Here $\alpha_\star$ is the rate found in
\cite{BBDGV, BDGV} for solutions emanating from integrable perturbations of
$U_{D,T}$. This fastest rate is the best constant in a Hardy-Poincar\'e
inequality (see \cite{BDGV}). This best constant is also the bottom
of the continuous spectrum of the linearization on a suitable weighted space
(see \cite{BBDGV, BDGV, DM}). In the figure below, the gray area indicates
the decay exponents that we obtain in
Theorem~\ref{mainthm} for $n=6$. The upper curve corresponds to the variational rate 
from \cite{BDGV}.

\vskip6mm 
\setbox2=\hbox{\epsfysize 5.5truecm\epsfbox[0 0 550 360]{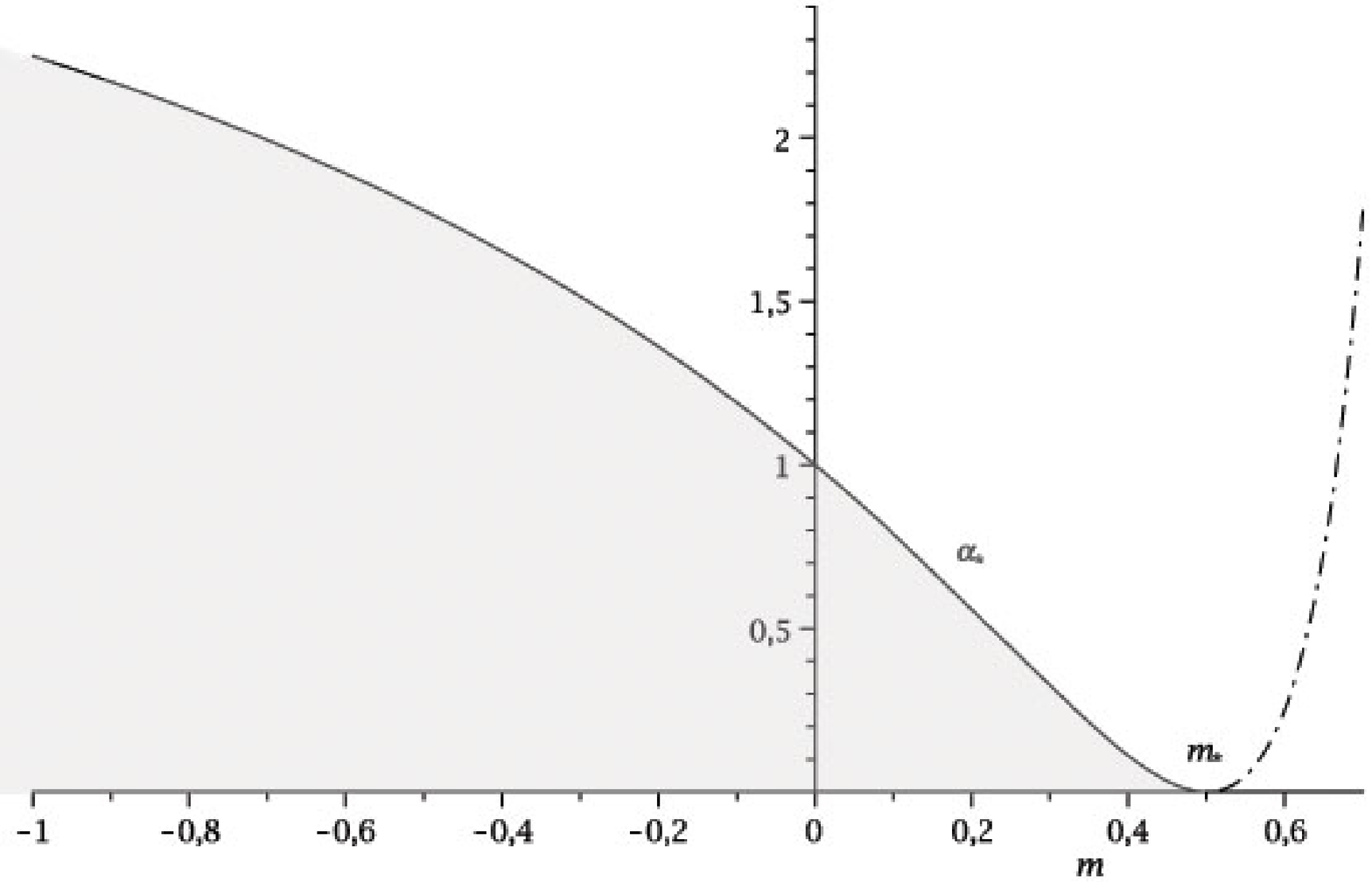}}
\vtop{
\centerline{\copy2}
\vskip3mm
\centerline{Figure 2.}
}
\vskip3mm plus5mm


\noindent (3)  
A continuum of rates of convergence to the
 singular Barenblatt profile $U_{0,T}$ was found in \cite{FVW}.
For the existence of continua of convergence rates to positive steady states, zero and
self-similar solutions of the Fujita equation
\[
u_t=\Delta u + u^p,\qquad x\in\R^n,
\]
see \cite{FWY1}, \cite{FWY2} and \cite{FWY3}, respectively.

\medskip

In Theorem~\ref{mainthm}, the assertion (i) is no longer true if
$l>l_\star$.
In fact, we have the following result about the   optimality of the range of
$l$.

\begin{theo}\label{subthm}
Let $m<m_*$, $n>2$. Assume that $D >0$ and
\[
0<v_0(x) < V_D(x),\qquad x\in\R^n
\]
or
\[
v_0(x)>V_D(x),\qquad x\in\R^n.
\]
Then for any $\eps>0$,
there exists $C_\eps>0$ such that the solution $v$ of {\rm (\ref{FPeqn})},
\eqref{eq:in}
satisfies
\be{1.2.1}
     \sup_{x\in\R^n} \Big | V_D(x)-v(x,t)\Big | \geq C_\eps
e^{-(\alpha_\star+\eps)t},
      \qquad t\ge 0.
\ee
\end{theo}

It follows from \eqref{1.2.1} that Theorem~2~(i) in \cite{BBDGV} is optimal
if $m<m_*$, $n>2$. The sharpness of the rate given by $\alpha_\star$
was discussed in \cite{BDGV} in terms of relative entropy which can be
written as
\[
{\cal F}[w]:= \frac{1}{1-m} \int_{\R^n}\left[w-1-\frac{1}{m}(w^m-1)\right]
V^m_D dx,\qquad w:=\frac{v}{V_D}\, .
\]
The statement on the sharp rate in \cite{BDGV} says that $\alpha =
\alpha_\star$ is the best possible rate for which
\[
{\cal F}[w(\cdot,t)]\leq {\cal F}[w(\cdot,0)]e^{-\alpha t}
\]
holds for all $t\geq 0$ if $V_{D_0}\leq v_0 \leq V_{D_1}$
for some $D_0>D>D_1>0$ and $v_0-V_D$ is integrable. Theorem~\ref{subthm}
implies that solutions starting from positive or negative perturbations
of $V_D$ cannot converge to $V_D$ (in $L^\infty$) at exponential rates
faster than $e^{-\alpha_\star t}$.

The paper is organised as follows. Section~2 contains some preliminaries,
in Section~3 we derive an upper bound for the rate of convergence
(Theorem~\ref{mainthm}~(i)) and in Section~4 we show that the rate we
found in Section~3 is optimal. We prove Theorem~\ref{subthm} in Section~5.
Section~6 is devoted to a few comments.

\mysection{Preliminaries}
Let us consider the initial value problem
 \be{19.1}
	\left\{ \begin{array}{l}
		(r^2+d)  \Big( \varphi_{rr} + \dfrac{n-1}{r} \varphi_{r}
\Big) - \mu r \varphi_{r}
		+ \alpha \varphi =0, \qquad r>0, \\[1mm]
		\varphi(0)=1, \quad \varphi_{r}(0)=0,
	\end{array} \right.
 \ee
with a parameter $d>0$, which is related to the linearized eigenvalue
problem at $V_D(x)$.
In order to study the behavior of the solution $\varphi=\varphi^d(r)$
of (\ref{19.1}), we first
consider properties of the  differential operator
\[
{\cal  L}W:=(r^2+d)  \Big( W_{rr} + \frac{n-1}{r} W_{r} \Big) - \mu r W_{r}
\, .
\]
We define
\be{l_alpha}
	l(\alpha):=\frac{n+\mu+2-\sqrt{(n-\mu-2)^2-4\alpha}}{2} >0
\ee
for $\alpha \in (0,\alpha_\star]$.
We note that for $\alpha \in (0,\alpha_\star)$, $l(\alpha)$ is the smaller
root of  the quadratic equation
\begin{equation} \label{ellalpha}
	  (l-\mu-2)(n-l)=\alpha.
\end{equation}

\begin{lem}\label{le:L}
Let $\alpha \in (0,\alpha_\star]$ and
$W^-(r):=(r^2+d)^{-k/2}$ with  $ k=l(\alpha)-\mu-2>0$.
Then  the inequality
\[
{\cal  L}W^- + \alpha W^- < 0
\]
holds for all $r>0$.
\end{lem}

\proof
Since
\[
\begin{aligned}
&W^-_r(r) = -k r(r^2+d)^{-k/2-1},
\\
& W^-_{rr}(r) = k(k+2) r^2(r^2+d)^{-k/2-2} - k(r^2+d)^{-k/2-1},
\end{aligned}
\]
we obtain
\[
\begin{aligned}
{\cal  L}W^-&=  k(k+2) r^2(r^2+d)^{-k/2-1} - k(r^2+d)^{-k/2}
  - (n-1) k (r^2+d)^{-k/2}
\\   &   \qquad  + \mu k r^2(r^2+d)^{-k/2-1}
   \\
   &= (r^2+d)^{-k/2-1} \Big (k(k+2) r^2- k(r^2+d) - (n-1)k(r^2+d)
    + \mu kr^2
   \Big )
   \\
   &= (r^2+d)^{-k/2-1} \Big ( (k^2-(n-2-\mu)k  ) r^2 - nkd
   \Big )
       \\
   &=   (k^2-(n-2-\mu)k ) W^-(r)
   -d (r^2+d)^{-k/2-1}   (k^2-(n-2-\mu)k   - nk  )  .
    \end{aligned}
\]
Here,  we have
\[
k^2-(n-2-\mu)k
=(l(\alpha)-\mu-2)(n-l(\alpha)) =-\alpha
\]
by $k=l(\alpha)-\mu-2$,
and
\[
k^2-(n-2-\mu)k   - nk   =
 l(\alpha)  (l(\alpha) -\mu-2) > 0
\]
by $l(\alpha)  > \mu + 2 $.
Hence  we obtain the desired inequality.
\qed

\

\begin{lem}\label{le:L2}
Let $\alpha \in (0,\alpha_\star)$ and
$W^+(r):=r^{-k} - r^{-j}$ with  $ k=l(\alpha)-\mu-2>0$ and
$j=l(\beta) - \mu - 2$.   If $\beta-\alpha>0$ is sufficiently small, then
there exists
$r_0>0$ such that   the inequality
\[
{\cal  L}W^+  + \alpha W^+ >0
\]
holds for $r>r_0$.
\end{lem}

\proof By direct computation, we have
\[
{\cal  L} r^{-k}  = - \alpha r^{-k} + O(r^{-k-2}),
\qquad {\cal  L}r^{-j}= - \beta r^{-j} + O(r^{-j-2}) \quad \mbox{ as } r \to
\infty.
\]
 Since $k<j<k+2$  if $\beta-\alpha>0$ is small,  we obtain
\[
{\cal  L}W^+ + \alpha W^+ =  (\beta-\alpha) r^{-j} + o(r^{-j})
\quad \mbox{ as } r \to \infty.
\]
This proves the lemma.
\qed

\

The following lemma shall be frequently used throughout the sequel.

\begin{lem}\label{lem19}
Let
$\alpha \in (0,\alpha_\star)$.
 Then for every $d>0$ the solution $\varphi=\varphi^d(r)$ of
\eqref{19.1}  is positive and decreasing in $r\in[0,\infty)$.
 Moreover, the solution has the following properties:
 \begin{description}
   \item[\rm(i)]  There exist positive constants $c$ and $C$ such that
 \be{19.2}
	c r^{-(l(\alpha)-\mu-2)} \le \varphi^d(r) \le C r^{-
(l(\alpha)-\mu-2)}, \qquad r \ge 1.
 \ee
\item[\rm(ii)]
$\varphi^d$ satisfies
\[
  \varphi^d_{rr}(r) + \frac{n-1}{r} \varphi^d_{r}(r) < 0, \qquad r>0.
\]
 \item[\rm(iii)]  $\varphi^d$ satisfies
\[
	0 \ge \frac{\varphi^d_{r}(r)}{\varphi^d(r)} \ge -
\big(l(\alpha)-\mu-2\big)
	\;  \frac{r}{r^2+d}, \qquad r>0.
\]
\end{description}
\end{lem}

\proof
Define $\rho(r):=r^{n-1}(r^2+d)^{-\mu/2}$.
Then the equation in (\ref{19.1}) is rewritten as
\begin{equation} \label{varphiaj}
 (\rho(r)\varphi^d_{r}(r))_r+ \alpha (r^2+d)^{-1} \rho(r) \varphi^d(r)=0,
\quad r>0.
\end{equation}
Similarly, it follows from  Lemma~\ref{le:L} that
$W^-(r):=(r^2+d)^{-k/2}$ with $ k=l(\alpha)-\mu-2>0$ satisfies
\begin{equation} \label{Waj}
(\rho(r)W^-_r(r))_r+ \alpha (r^2+d)^{-1} \rho(r) W^-(r)<0, \quad r>0.
\end{equation}
Multiplying (\ref{varphiaj}) by $W^-(r)$ and (\ref{Waj}) by $\varphi^d(r)$,
taking the
difference and integrating it by parts on $[0,r]$, we obtain
\begin{equation} \label{green}
 \rho(r) \big\{ \varphi^d_r(r) W^-(r) - \varphi^d(r) W^-_r(r) \big\} > 0
\end{equation}
as long as $\varphi^d(r)>0$.
This implies that   $\varphi^d(r)/W^-(r)$ is increasing
in $r$,  and hence $\varphi^d(r)>0$ for all $r>0$.
Hence there exist $C_1>0$ such that
$C_1 W^-(r) < \varphi^d(r)
$
for all $r>0$.

Similarly, from Lemma~\ref{le:L2} and  (\ref{varphiaj}), we have
\begin{equation} \label{Wplus}
\Big[ \rho( \varphi^d_r  W^+ - \varphi^d W^+_r ) \Big]_{r_0}^r
 < -  \alpha \int_{r_0}^r (s^2+d)^{-1} \rho(s) \varphi^d(s) W^+(s)  ds.
\end{equation}
Here the integrand satisfies
\[
\begin{aligned}
(s^2+d)^{-1} \rho(s) \varphi^d(s) W^+(s)
& \geq C_1 (s^2+d)^{-1} \rho(s) W^-(s) W^+(s)  \\
& = (C_1 + o(1))  s^{-2}  s^{n-1-\mu} s^{-k} s^{-k}  \\
&= (C_1 + o(1)) s^{n+\mu+1 - 2l(\alpha)} \quad \mbox{ as } s\to\infty.
\end{aligned}
\]
Since $n+\mu+1 - 2l(\alpha)>-1$, the right-hand side of (\ref{Wplus})
tends to $-\infty$.
This implies that  $\varphi^d(r)/W^+(r)$ is decreasing for large $r$, and
hence there exists $C_2>0$ such that
$  \varphi^d(r) < C_2 W^+(r) $
for large $r$.  This completes the proof of (i).

Integrating (\ref{varphiaj}) on $[0,r]$ and using the initial condition, we
obtain
\[
\rho(r)\varphi^d_{r}(r) = - \alpha \dint_0^r (s^2+d)^{-1} \rho(s)
\varphi^d(s) ds<0,
\]
which implies $\varphi^d_{r}(r)<0$ for all $r>0$,  Hence by
 the equation in (\ref{19.1}),
(ii)  is proved.

Finally, by (\ref{green}), we have
\[
 \dfrac{\varphi^d_r(r)}{\varphi^d(r)}  \geq \dfrac{W^-_r(r)}{W^-(r)}
  = - \big(l(\alpha)-\mu-2\big)
	\;  \frac{r}{r^2+d}.
\]
Thus (iii) is proved.
\qed

\

\begin{lem}\label{le:alphazero}
If $\alpha>\alpha_\star$, then the solution $\varphi=\varphi^d(r)$
of \eqref{19.1} changes sign.
 \end{lem}

\proof
We will derive a contradiction by assuming  that $\varphi^d(r)>0$ for all
$r>0$.
Set $W_\star(r):=r^{-k}$ with  $ k=l(\alpha_\star)-\mu-2=(n-\mu-2)/2$.
Then   we have
\[
{\cal  L} W_\star = -\alpha_\star r^{-k} + O(r^{-k-2}) \quad \mbox{ as } r
\to \infty.
\]
This implies that for any fixed $a \in (0,\alpha-\alpha_\star)$,
there exists $r_1>0$  such that
\[
(\rho (r) (W_\star(r))_r)_r  + \alpha  (r^2+d)^{-1} \rho(r) W_\star(r) >
   a  (r^2+d)^{-1} \rho(r) W_\star(r),   \qquad  r\geq r_1.
\]
 Then from (\ref{varphiaj}), we  obtain
\begin{equation} \label{eq:wphi}
 \Big [ \rho (  \varphi^d_r W_\star - \varphi^d (W_\star)_r  )  \Big
]_{r_1}^r
< - a   \dint_{r_1}^r  (s^2+d)^{-1} \rho(s)\varphi^d(s)  W_\star(s)  ds.
\end{equation}

Suppose here that $(\varphi^d(r)/W_\star(r))_r>0$ holds for $r \geq r_1$.
Then there exists a constant $c>0$ such that $\varphi_d(r)>cr^{-k}$
for $r \geq r_1$.  In this case, the integrand in (\ref{eq:wphi}) satisfies
\[
(s^2+d)^{-1} \rho(s)\varphi^d(s)  W_\star(s) \geq  (c+o(1)) s^{n-\mu-2}
s^{-2k}
 = (c+o(1))s^{-1} \quad \mbox{ as } s \to \infty.
\]
Hence the right-hand side of  (\ref{eq:wphi}) tends to $-\infty$ as $r \to
\infty$,
which contradicts the assumption
that  $(\varphi^d(r)/W_\star(r))_r>0$ for $r \geq r_1$.

Thus we may assume without loss of generality
that $(\varphi^d(r)/W_\star(r))_r \leq 0$ at $r = r_1$.
Then by (\ref{eq:wphi}), we have
\begin{equation} \label{pwgr}
 \rho(r)  \big \{  \varphi^d_r (r) W_\star(r) - \varphi^d (r)(W_\star(r))_r
\big  \}
 <-  a   \dint_{r_1}^r  (s^2+d)^{-1} \rho(s)\varphi^d(s)  W_\star(s)  ds
\end{equation}
for $r \geq r_1$.
Hence  $(\varphi^d(r)/W_\star(r))_r <0$ for  $r \geq r_1$
so that $\varphi^d(r)/W_\star(r)$ converges to a
nonnegative  constant. This implies that
there exists a sequence $\{ r_i \}$ such that $r_i \to \infty$ as $i \to
\infty$
and
$(\varphi^d(r)/W_\star(r))_r \to 0$ along the sequence.  Then
the left-hand side of (\ref{pwgr}) tends to 0
along the sequence, since $\rho(r)$ and $W_\ast(r)^{-2}$ are of the same
order
as $r \to \infty$, whereas the right-hand of
is negative and decreasing in $r>r_1$.
This contradiction shows that $\varphi^d(r)$ must change sign at some
finite $r$.   \qed

\

Let $\parab$ be the differential operator
\be{2.100}
	\parab w := w_t- \Big( (w^{m-1}w_r)_{r} + \frac{n-1}{r} w^{m-1}w_r \Big)
	- \mu r w_r - \mu n w , \qquad r>0, \ t>0.
\ee
In the following lemma we evaluate $\parab w$
for a suitable auxiliary function $w$ which will be used frequently in the
sequel.

\begin{lem}\label{lem18}
 Let $D>0$, and   $y:  [0,\infty) \to \R$ and $\psi: \, (0,\infty) \to
(0,\infty)$ be
 any  smooth functions.
 Then the function
 \bas
	w(r,t):=\Big(r^2+D+y(t)  \psi(r) \Big)^{- \mu/2}, \qquad r>0, \ t>0,
 \eas
 satisfies
 \be{18.1}
	\parab w= \frac{\mu}{2}  y(t)  \Big(r^2+D+y(t)
\psi(r)\Big)^{-(\mu+2)/2}
	 \AD[y(t)] \psi, \qquad r>0, \ t>0,
 \ee
 where
 \bea{18.2}
	\AD[y(t)] \psi &:=& (r^2+D) \,\Big(\psi_{rr} + \frac{n-1}{r} \psi_r
\Big)
	-\mu r \psi_r - \frac{y'(t)}{y(t)} \ \psi \nn\\
	& & - y(t) \, \bigg\{ -\psi \, \Big(\psi_{rr} + \frac{n-1}{r}
\psi_r\Big)
	+ \frac{\mu}{2} (\psi_r)^2 \bigg\}
 \eea
 for $r>0$ and $t>0$.
\end{lem}

\proof
 We compute
 \bas
	w_t=-\frac{\mu}{2} \Big(r^2+D+y(t) \ \psi(r) \Big)^{-(\mu+2)/2}  y'
\psi
 \eas
 and
 \bas
	w_r = - \frac{\mu}{2}  \Big(r^2+D+y(t)  \psi(r) \Big)^{-(\mu+2)/2}
(2r+y\psi_r),
 \eas
 and using the relation $m=(\mu-2)/\mu$ we find that
\[
\begin{aligned}
	w^{m-1}w_r &= \frac{\mu}{\mu-2} \Big( (r^2+D+y
\psi)^{-(\mu-2)/2}\Big)_r \\
	&= -\frac{\mu}{2} \Big(r^2+D+y(t)  \psi(r) \Big)^{-\mu/2}  (2r+y
\psi_r)
\end{aligned}
\]
 as well as
\[
\begin{aligned}
	(w^{m-1}w_r)_r &= \frac{\mu^2}{4}
	\Big(r^2+D+y(t)  \psi(r) \Big)^{-(\mu+2)/2}  (2r+y\psi_r)^2 \\
	&  -\frac{\mu}{2} \Big(r^2+D+y(t)  \psi(r) \Big)^{-\mu/2}
(2+y\psi_{rr})
\end{aligned}
\]
 for $r>0$ and $t>0$. Therefore,
\[
\begin{aligned}
	\parab w  = \frac{\mu}{2} \Big(r^2+D+y(t)  \psi(r)
\Big)^{-(\mu+2)/2}  \Big\{ &
		-y' \psi 	- \frac{\mu}{2} (2r+y\psi_r)^2 \\
		& + (r^2+D+y\psi)  (2+y\psi_{rr}) \\
	&
	+ \frac{n-1}{r} (r^2+D+y\psi)  (2r+y\psi_r) \\
	&
	+ \mu r (2r+y\psi_r) 	-2n (r^2+D+y\psi) \Big\} \\[2mm]
	= \frac{\mu}{2} \Big(r^2+D+y(t)  \psi(r) \Big)^{-(\mu+2)/2}
\Big\{ & - y'\psi
		-2\mu r^2 - 2\mu yr\psi_r - \frac{\mu}{2} y^2 (\psi_r)^2
		\\
		& + 2(r^2+D+y\psi)  	+ y(r^2+d+y\psi) \psi_{rr} \\
	&
	+2(n-1) (r^2+D+y\psi) \\
	&  + y (r^2+D+y\psi)  \frac{n-1}{r} \psi_r 	+2\mu r^2 + \mu
yr\psi_r \\
	&
	-2n(r^2+D+y\psi) \Big\},
\end{aligned}
\]
 which after obvious simplifications yields (\ref{18.1}).
\qed

\mysection{Upper bound for the convergence rate}
The goal of this section is to provide an upper estimate for
the distance between the solution $v$ of (\ref{FPeqn})
and $V_D$ under the assumption that $|v_0(r)-V_D(r)| \le c r^{-l}$
for all $r\ge 1$ and some $c>0$ and $l \in (\mu+2,l_\star)$ with
$l_\star$ given in Theorem~\ref{mainthm}.

\subsection{Estimating $v(\cdot,t)-V_D$ from above}
We first consider the {\em positive part} of the  distance,
that is, we seek an upper bound for
$\sup_{r\ge 0} (v(r,t)-V_D(r))$ for large $t$. This will be accomplished
along
a three-step procedure: In a first step, in Lemma~\ref{lem27}
we shall establish that this term will be sufficiently small, but at a rate
possibly
smaller than the desired one. Next, in Lemma~\ref{lem30} we shall
make sure that the spatial decay rate of the `positive tail' $
(v(\cdot,t)-V_D)_+$ is essentially
preserved during the evolution. Finally, the proof of Lemma~\ref{lem31}
will be based upon the observation that if at a certain
time the solution is close enough to $V_D$ {\em and} its positive tail
decays sufficiently fast in space, then $(v(\cdot,t)-V_D)_+$ will
decay in time at the desired rate.

In order to avoid redundant arguments, we start with the following
lemma that will be used twice in Lemma~\ref{lem27} and
in Lemma~\ref{lem30}.

\begin{lem}\label{lem36}
Let $D>0, \delta \in (0,D), B>0$ and $\alpha \in (0,\alpha_\star)$
with $\alpha_\star$ given by {\rm (\ref{alpha_star})}, and
suppose that $\varphi^\delta$ solves {\rm (\ref{19.1})} with $d=\delta$.
Moreover, assume that $\eta$ is a real number satisfying
\be{36.2}
	\eta \le \alpha - \frac{\mu(D-\delta)}{2\delta}
\Big(l(\alpha)-\mu-2\Big)^2,
\ee
where $l(\alpha)$ is defined by \eqref{l_alpha}.
Then the function $\ov$ defined by
\be{36.444}
	\ov(r,t) := \min \bigg\{ V_\delta(r) \, , \, \Big(r^2+D-Be^{-\eta t}
	\varphi^\delta(r)\Big)_+^{- \mu/2} \bigg\},
	\qquad r\ge 0, \ t\ge 0,
\ee
satisfies $\parab \ov \ge 0$ in the Nagumo sense in $(0,\infty) \times
(0,\infty)$.
\end{lem}

\proof
Since $V_\delta$ is a solution of (\ref{FPeqn}),
we only need to make sure that $\parab \ov \ge 0$ is valid at each point
from the set
\bas
	S:=\Big\{ (r,t) \in (0,\infty) \times (0,\infty) \ \Big| \ \ov(r,t)
< V_\delta(r) \Big\}.
\eas
In view of Lemma~\ref{lem18} this amounts to showing that at such points
$\AD[-Be^{-\eta t}]  \varphi^\delta \le 0$ holds, where
\bea{36.44}
	\AD[-Be^{-\eta t}]  \varphi^\delta &=& (r^2+D)  \Big(
\varphi^\delta_{rr}
	+ \frac{n-1}{r} \varphi^\delta_r \Big)
	- \mu r \varphi^\delta_r + \eta \varphi^\delta \nn\\
	& & + Be^{-\eta t}  \bigg\{ - \varphi^\delta
\Big(\varphi^\delta_{rr} + \frac{n-1}{r}
	 \varphi^\delta_r \Big)
	+\frac{\mu}{2} (\varphi^\delta_r)^2 \bigg\}.
\eea
To see this, we first note that the nonlinear terms
$- \varphi^\delta  \big(\varphi^\delta_{rr} + \frac{n-1}{r}
\varphi^\delta_r \big)$
and $(\varphi^{\delta}_{r})^2$ in (\ref{36.44})
are both positive on $(0,\infty)$ by Lemma~\ref{lem19}. In order to
estimate
 them appropriately, we make use of the fact that whenever
$\ov(r,t) < V_\delta(r)$, due to (\ref{36.444}) we know that
\[
	r^2+\delta = V_\delta(r) ^{-2/\mu}\\
	< \ov(r,t) ^{-2/\mu}\\
	= r^2+D-B e^{-\eta t}  \varphi^\delta(r)
\]
and hence
\bas
	B e^{-\eta t} < \frac{D-\delta}{\varphi^\delta(r)}.
\eas
Using this in (\ref{36.44}) yields
\[
\begin{aligned}
	\AD[-Be^{-\eta t}]
	 \varphi^\delta &\le (r^2+D)  \Big( \varphi^\delta_{ rr}
	 + \frac{n-1}{r} \varphi^\delta_r \Big)
	- \mu r \varphi^\delta_r+ \eta \varphi^\delta \\
	&  - (D-\delta)  \Big( \varphi^\delta_{rr} + \frac{n-1}{r}
\varphi^\delta_r\Big)
	+\frac{\mu(D-\delta)(\varphi^\delta_r)^2}{2\varphi^\delta}    \\
	&= (r^2+\delta)  \Big( \varphi^\delta_{rr} + \frac{n-1}{r}
\varphi^\delta_r\Big)
		-\mu r \varphi^\delta_r+ \eta \varphi^\delta
		+ \frac{\mu(D-\delta)(\varphi^\delta_r)^2}{2\varphi^\delta}
\end{aligned}
\]
at such points.
According to the ODE in (\ref{19.1}), we thus obtain that
\bea{36.6}
	\AD[-Be^{-\eta t}]  \varphi^\delta &\le& (-\alpha \varphi^\delta +
\mu r \varphi^\delta_r)
	-\mu r \varphi^\delta_r+ \eta \varphi^\delta
	+ \frac{\mu(D-\delta)(\varphi^\delta_r)^2}{2\varphi^\delta}    \nn\\
	&=& -(\alpha-\eta) \varphi^\delta
	+ \frac{\mu(D-\delta)(\varphi^\delta_r)^2}{2\varphi^\delta}
\eea
is satisfied at all points in $S$.
Now from Lemma~\ref{lem19} (iii) we know that
\[
\begin{aligned}
	\left ( \frac{\varphi^\delta_r}{\varphi^\delta} \right )^2&
	\le \Big(l(\alpha)-\mu-2\Big)^2  \frac{r^2}{(r^2+\delta)^2} \\
	&\le \Big(l(\alpha)-\mu-2\Big)^2  \frac{1}{r^2+\delta} \\
	&\le \frac{1}{\delta} \Big(l(\alpha)-\mu-2\Big)^2   \qquad \mbox{for
all } r>0,
\end{aligned}
\]
which inserted into (\ref{36.6}) yields
\[
	\AD[-Be^{-\eta t}]
	 \varphi^\delta  \le  \bigg\{ -(\alpha-\eta) +
\frac{\mu(D-\delta)}{2\delta}
	  \Big(l(\alpha)-\mu-2\Big)^2
	  \bigg\}
	 \varphi^\delta
	\qquad \mbox{in } S.
\]
This immediately proves that under the assumption
 (\ref{36.2}), $\AD[-Be^{-\eta t}]  \varphi^\delta \le 0$ holds throughout
$S$.
\qed

\

As announced, our three-step procedure is launched by
deriving an estimate from above.

\begin{lem}\label{lem27}
Suppose that there exist $\delta>0$, $D>\delta$, $c>0$ and   $l \in (
\mu+2,l_\star)$ such that
\be{h0}
	0 < v_0(r) \le V_\delta(r) \qquad \mbox{for all } r \ge 0
\ee
and
\be{27.1}
	v_0(r) \le V_D(r) + cr^{-l} \qquad \mbox{for all } r\ge 1.
\ee
Then there exist $C>0, \eta>0$ and $l_0 \in (\mu+2,l]$
such that the solution $v$ of (\ref{FPeqn}) satisfies
\be{27.2}
	v(r,t) \le V_D(r) + C \, e^{-\eta t}  (r+1)^{-l_0} \qquad \mbox{for
all $r\ge 0$ and } t\ge 0.
\ee
\end{lem}

\proof
Since our assumption $m<(n-4)/(n-2)$ implies $\mu+2<n$,
 we can pick some $l_0>\mu+2$ close enough to $\mu+2$ such that
$l_0 \le l$, $l_0<l_\star$ and
\be{27.3}
	\frac{l_0-\mu-2}{n-l_0} \le \frac{\delta}{\mu(D-\delta)}.
\ee
Then $\alpha:=(l_0-\mu-2)(n-l_0)$ is positive, so that we can pick some
$\eta \in(0, \alpha/2)$, and by (\ref{l_alpha}) $l_0$ is given by
\be{27.4}
	l_0=l(\alpha).
\ee
Now from the hypothesis (\ref{27.1}) and the convexity
 of $0< z \mapsto z^{-2/\mu}$ we obtain
\[
\begin{aligned}
	v_0(r)^{-2/\mu} &\ge \Big(V_D(r)+cr^{-l}\Big)^{-2/\mu} \\
	&\ge V_D(r)^{-2/\mu}  \Big(1-\Big(\frac{2c}{\mu}\Big)
\frac{r^{-l}}{V_D(r)} \Big) \\
	&= r^2+D - \Big(\frac{2c}{\mu}\Big)  (r^2+D)^{ (\mu+2)/2}  r^{-l}
	\qquad \mbox{for all } r\ge 1,
\end{aligned}
\]
which combined with the positivity of $V_D$ and (\ref{19.2}) easily yields
that
\be{28.1}
	V_D^{-2/\mu}(r) \ge r^2+D-B \varphi^\delta(r) \qquad \mbox{for all }
r\ge 0
\ee
holds with some suitably large $B>0$, where $\varphi^\delta$
denotes the solution of (\ref{19.1}) with $d=\delta$.
In view of (\ref{h0}), this implies that the function $\ov$ defined by
(\ref{36.444}) dominates $v$ at $t=0$. Since moreover
\[
	\frac{\mu(D-\delta)}{2\delta}  \Big(l(\alpha)-\mu-2\Big)^2
	=\frac{\mu(D-\delta)}{2\delta}  \Big(l_0-\mu-2\Big)^2
	=\frac{\alpha\mu(D-\delta)(l_0-\mu-2)}{2(n-l_0)\delta} \le
\frac{\alpha}{2}
\]
by (\ref{27.4}) and (\ref{27.3}), Lemma~\ref{lem36} yields that $\ov$
 is a subsolution of (\ref{FPeqn}) and hence we have
$v \le \ov$ by comparison. Now (\ref{27.2}) is a straightforward
 consequence of the definition of $\ov$:
We let $c_1>0$ be such that $(1-z)^{-\mu/2} \le 1+c_1 z$
for all $z \in [0,B/D]$ and observe that
since $\varphi^\delta_r\le 0$ we have
\[
 \frac{B e^{-\eta t}  \varphi^\delta(r)}{r^2+D} \le \frac{B}{D}
 \]
  for all $r\ge 0$
and $t\ge 0$. Hence,
\bea{28.2}
	v(r,t) &\le& \ov(r,t) \nn\\
	&\le& \Big( r^2+D-Be^{-\eta t}  \varphi^\delta(r)\Big)^{-\mu/2}
\nn\\
	&=& (r^2+D)^{- \mu/2}
		\bigg( 1 - \frac{B e^{-\eta t}  \varphi^\delta(r)}{r^2+D}
\bigg)^{-\mu/2} \nn\\
	&\le& (r^2+D)^{-\mu/2} + c_1 B e^{-\eta t}  (r^2+D)^{-(\mu+2)/2}
\varphi^\delta(r)
\eea
for all $r\ge 0$ and $t\ge 0$. By means of the second
 inequality in (\ref{19.2}) we easily arrive at (\ref{27.2}).
\qed

\

We next assert that the decay of $(v(\cdot,t)-V_D)_+$ does not
change in time significantly.

\begin{lem}\label{lem30}
Under the assumption of Lemma~\ref{lem27},
there exist $C>0$ and $\kappa>0$
such that the solution $v$ of {\rm (\ref{FPeqn})} satisfies
\be{30.2}
	v(r,t) \le V_D(r) + C e^{\kappa t}  (r+1)^{-l} \qquad \mbox{for all
$r\ge 0$ and } t\ge 0.
\ee
\end{lem}

\proof
Since $l \in (\mu+2,l_\star)$, the number $\alpha:=(l-\mu-2)(n-l)$
 is positive and satisfies $\alpha<\alpha_\star$
with $\alpha_\star$ as in (\ref{alpha_star}), and moreover the
number $l(\alpha)$ in (\ref{l_alpha}) coincides with $l$.
We now choose positive constants $\kappa$ and $B$ such that
\be{30.3}
	\kappa \ge \frac{\mu(D-\delta)(l-\mu-2)^2}{2\delta}-\alpha
\ee
and
\be{30.4}
	v_0^{-2/\mu}(r) \ge r^2+D-B\varphi^\delta(r) \qquad \mbox{for all }
r\ge 0,
\ee
where the latter is possible because of (\ref{27.1}) and (\ref{19.2})
(cf.~the argument leading to (\ref{28.1})).
Thus, defining $\ov$ by (\ref{36.444}) with $\eta:=-\kappa$, (\ref{30.4})
states that $v_0(r) \le \ov(r,0)$ for all $r\ge 0$,
whereas (\ref{30.3}) ensures that (\ref{36.2}) holds.
Therefore, Lemma~\ref{lem36} in conjunction with the comparison principle
guarantees that $v \le \ov$ in $[0,\infty) \times [0,\infty)$,
whence (\ref{30.2}) easily results from the definition of $\ov$ (see
(\ref{28.1})).
\qed

\

We are now in the position to prove the main result of this section.

\begin{lem}\label{lem31}
Under the assumption of Lemma~\ref{lem27},
the solution $v$ of {\rm (\ref{FPeqn})} satisfies
\be{31.2}
	v(r,t) \le V_D(r) + C e^{-(l-\mu-2)(n-l)t}  (r+1)^{-(l-\mu-2)}
	\qquad \mbox{for all $r\ge 0$ and } t\ge 0
\ee
with some $C>0$. In particular, $v$ satisfies
\be{31.3}
	\sup_{r\ge 0} \Big(v(r,t)-V_D(r)\Big) \le C e^{-(l-\mu-2)(n-l)t}
\qquad \mbox{for all } t\ge 0.
\ee
\end{lem}

\proof
Again, $\alpha:=(l-\mu-2)(n-l)$ belongs to $(0,\alpha_\star)$ and is such
that $l(\alpha)=l$.
Now taking any small $\delta>0$ such that
\be{31.4}
	\delta < \frac{(n+\mu-l)D}{n+\mu-l+\frac{\mu}{2}(l-\mu-2)},
\ee
from Lemma~\ref{lem27} we infer the existence of some large $t_0>0$ such
that
\be{31.5}
	v(r,t_0) \le V_{D-\delta}(r) \qquad \mbox{for all } r \ge 0,
\ee
because the difference $ V_{D-\delta}(r)  -  V_D(r) $ decays with the order
$r^{-(\mu+2)}$.
Then in view of Lemma~\ref{lem30} it is possible to fix a large number
$c_1>0$ such that
\be{31.6}
	v(r,t_0) \le V_D(r) + c_1 r^{-l} \qquad \mbox{for all } r \ge 1.
\ee
We finally choose $B$ appropriately large satisfying
\be{31.7}
	B \ge e^{\alpha t_0}  (D+1)
\ee
and
\be{31.8}
	B \ge \frac{2c_1}{\mu} e^{\alpha t_0}  (D+1)^{(\mu+2)/2}
\ee
and let
\be{31.9}
	\ov(r,t) := \min \bigg\{ V_{D-\delta}(r) \, , \,
	\Big( r^2+D-Be^{-\alpha t} r^{-(l-\mu-2)} \Big)^{-\mu/2}_+ \bigg\},
	\qquad r\ge 0, \ t\ge 0.
\ee
We first claim that $\ov \ge v$ holds for all $r\ge 0$ at $t=t_0$.
Indeed, if $r\le 1$ then
\bas
	r^2+D-Be^{-\alpha t_0}  r^{-(l-\mu-2)} \le 1+D-B e^{-\alpha t_0} \le
0
\eas
by (\ref{31.7}), hence $\ov(r,t_0) = V_{D-\delta}(r) \ge v(r,t_0)$ by
(\ref{31.5}).
On the other hand, if $r>1$ then by convexity of $0<z \mapsto z^{-2/\mu}$
and (\ref{31.6}),
\[
\begin{aligned}
	v(r,t_0)^{-2/\mu} &\ge \Big(V_D(r)+c r^{-l}\Big)^{-2/\mu} \\
	&= V_D(r)^{-2/\mu}  \Big( 1+\frac{cr^{-l}}{V_D(r)}\Big)^{-2/\mu} \\
	&\ge V_D(r)^{-2/\mu}  - \frac{2c}{\mu} V_D(r)^{-2/\mu -1} r^{-l} \\
	&= r^2 + D - \frac{2c}{\mu}  \Big(\frac{r^2+D}{r^2}
\Big)^{(\mu+2)/2}  r^{-(l-\mu-2)} \\
	&\ge r^2+D-\frac{2c}{\mu}  (D+1)^{(\mu+2)/2} r^{-(l-\mu-2)} \\
	&\ge r^2+D-B e^{-\alpha t_0}  r^{-(l-\mu-2)}
	\qquad \mbox{for all } r>1
\end{aligned}
\]
in view of (\ref{31.8}).

Having thereby asserted the desired ordering at $t=t_0$,
let us next make sure that $\ov$ is a supersolution of the
PDE in (\ref{FPeqn}) for $r>0$ and $t>t_0$. To achieve this,
since $V_{D-\delta}$ solves (\ref{FPeqn}), we only 
need to check that $\AD[-Be^{-\alpha t}]  \psi \le 0$ in $S$ is valid with
$\psi(r):=r^{-(l-\mu-2)}$ and
$S:=\{(r,t) \in (0,\infty) \times (t_0,\infty) \ | \ \ov(r,t) <
V_{D-\delta}(r) \}$, where
\[
\begin{aligned}
	\AD [-Be^{-\alpha t}]
	\psi&= (r^2+D)  \Big( \psi_{rr} + \frac{n-1}{r} \psi_r \Big) - \mu r
\psi_r + \alpha \psi \\
	& \qquad
	 + B e^{-\alpha t}  \Big\{ - \psi  \Big( \psi_{rr} + \frac{n-1}{r}
\psi_r \Big)
	+ \frac{\mu}{2} \psi_r^2 \Big\}.
\end{aligned}
\]
Differentiating $\psi$, we immediately obtain
\[
\begin{aligned}
	\AD[-Be^{-\alpha t}] \psi &= - (r^2+D) (l-\mu-2)(n+\mu-l)
r^{-(l-\mu)} \\
	& \qquad  + \mu(l-\mu-2) r^{-(l-\mu-2)} + \alpha r^{-(l-\mu-2)} \\
	& \qquad + B e^{-\alpha t}  \Big\{
	r^{-(l-\mu-2)}  (l-\mu-2)(n+\mu-l) r^{-(l-\mu)} \\
	&  \hspace*{20mm}
	\qquad+ \frac{\mu}{2} \Big( (l-\mu-2)^2 r^{-(l-\mu-1)}\Big)^2 \Big\}
\\[1mm]
	&= - (l-\mu-2)(n+\mu-l)  (r^2+D)  r^{-(l-\mu)} \\
	&  \qquad+ \Big( \mu(l-\mu-2) + \alpha \Big)  r^{-(l-\mu-2)} \\
	& \qquad+ B e^{-\alpha t}  (l-\mu-2) \Big( n+\mu-l +
\frac{\mu}{2}(l-\mu-2) \Big)  r^{-2(l-\mu-1)},
\end{aligned}
\]
so that by definition of $\alpha$ we see that
\[
\begin{aligned}
	\AD[-Be^{-\alpha t}] \psi &= -(l-\mu-2)(n+\mu-l)  D r^{-(l-\mu)} \\
	& \qquad + B e^{-\alpha t}  (l-\mu-2) \Big( n+\mu-l +
\frac{\mu}{2}(l-\mu-2) \Big)  r^{-2(l-\mu-1)}
\end{aligned}
\]
for $r>0$ and $t>0$. Now if $(r,t) \in S$ then, by (\ref{31.9}),
\bas
	B e^{-\alpha t} < \delta r^{l-\mu-2},
\eas
and therefore at these points we have
\[
\begin{aligned}
	\AD[-Be^{-\alpha t}] \psi &\le -(l-\mu-2)(n+\mu-l)  D r^{-(l-\mu)}
\\
	&  \qquad +\delta(l-\mu-2) \Big( n+\mu-l
	+ \frac{\mu}{2}(l-\mu-2) \Big)  r^{-(l-\mu)} \\
	&= (l-\mu-2) \Big\{ -(n+\mu-l)D \\
& \qquad + \delta \Big( n+\mu-l
	+ \frac{\mu}{2}(l-\mu-2) \Big) \Big\}  r^{-(l-\mu)} \\
	&\le 0
\end{aligned}
\]
because of our smallness condition (\ref{31.4}) on $\delta$.
Consequently, by comparison we conclude that $v \le \ov$
 for $r\ge 0$ and $t\ge t_0$, and hence (\ref{31.2}) and
(\ref{31.3}) easily result from the definition (\ref{31.9}) of $\ov$.
\qed

\subsection{Estimating $v(\cdot,t)-V_D$ from below}
A corresponding statement on convergence {\em from below}
can be obtained much more directly by a one-step comparison
argument.

\begin{lem}\label{lem33}
Suppose that $v_0>0$ on $[0,\infty)$, and that there
exist $D>0$, $c>0$ and $l\in (\mu+2,l_\star)$ such that
\be{33.1}
	v_0(r) \ge V_D(r) - c r^{-l} \qquad \mbox{for all } r>0.
\ee
Then there exists $C>0$ such that the solution $v$ of \eqref{FPeqn}
satisfies
\be{33.2}
	v(r,t) \ge V_D(r) - C e^{-(l-\mu-2)(n-l)t}  (r+1)^{-l}
	\qquad \mbox{for all $r\ge 0$ and } t\ge 0
\ee
and thus
\bas
	\sup_{r\ge 0} \Big(V_D(r)-v(r,t) \Big) \le C e^{-(l-\mu-2)(n-l)t}
	\qquad \mbox{for all } t\ge 0.
\eas
\end{lem}

\proof
Since $l>\mu+2>\mu$, it is possible to fix $r_0 \ge 1$
such that with $c$ as in (\ref{33.1}) we have
\be{33.3}
	c r^{-l} \le \frac{1}{2} V_D(r) = \frac{1}{2} (r^2+D)^{-\mu/2}
\qquad \mbox{for all } r>r_0.
\ee
Then according to Lemma~\ref{lem19} there exists $c_1>0$
 such that the solution $\varphi^D$ of (\ref{19.1}) with $d=D$ satisfies
\be{33.4}
	\varphi^D(r) \ge c_1 r^{-(l-\mu-2)} \qquad \mbox{for all  } r>r_0,
\ee
and using that $v_0>0$ on $[0,\infty)$ we can pick $c_2>0$ such that
\be{33.5}
	v_0(r) \ge c_2 \qquad \mbox{for all } r\in [0,r_0].
\ee
We next let $c_3>0$ be such that
\be{33.6}
	(1-z)^{-\frac{2}{\mu}} \le 1 + c_3 z \qquad \mbox{for all } z\in
[0,1/2]
\ee
and finally fix $B>0$ large enough satisfying
 \be{33.7}
	B \ge \Big ({c_2^{2/\mu} \varphi^D(r_0)}\Big)^{-1}
\ee
as well as
\be{33.8}
	B \ge \frac{\gamma (D+1)^{(\mu+2)/2}}{c_1}, \qquad \gamma:=cc_3.
\ee
Then, again writing $\alpha:=(l-\mu-2)(n-l)>0$, from Lemma~\ref{lem18} we
know that
\bas
	\uv(r,t):= \Big( r^2+D+B e^{-\alpha t}  \varphi^D(r) \Big)^{-\mu/2}
,
	\qquad r\ge 0, \ t\ge 0,
\eas
defines a subsolution of (\ref{FPeqn}) if and only
if $\AD[Be^{-\alpha t}] \varphi^D \le 0$, where
\[
\begin{aligned}
	\AD[Be^{-\alpha t}] \varphi^D &= (r^2+D)
	\Big( \varphi^D_{rr} + \frac{n-1}{r} \varphi^D_{ r} \Big)
	-\mu r \varphi^D_{r} + \alpha \varphi^D \\
	   & \qquad
	   - B e^{-\alpha t}  \bigg\{ - \varphi^D
	   \Big( \varphi^D_{rr} + \frac{n-1}{r} \varphi^D_{ r} \Big)
	+\frac{\mu}{2} (\varphi^D_{ r})^2 \bigg\}, \quad r>0, \ t>0.
\end{aligned}
\]
Recalling that Lemma~\ref{lem19} yields
$\varphi^D_{ rr} + \frac{n-1}{r} \varphi^D_{ r} \le 0$ on $(0,\infty)$,
in view of (\ref{18.1}) we directly see that indeed
\[
	\AD [Be^{-\alpha t}]\varphi^D \le (r^2+D)
	\Big( \varphi^D_{ rr} + \frac{n-1}{r} \varphi^D_{ r} \Big)
	-\mu r \varphi^D_{ r} + \alpha \varphi^D =0
\]
for all $r>0$ and $t>0$.
In order to use $\uv$ as a comparison function, we next claim that
\be{33.10}
	\uv(r,0) \le v_0(r) \qquad \mbox{for all } r\ge 0.
\ee
In fact, for small $r$, (\ref{33.5}) and (\ref{33.7}) tell us that
\[
\begin{aligned}
	\frac{\uv(r,0)}{v_0(r)} &\le \frac{\uv(r,0)}{c_2}
	= \frac{1}{c_2}  \Big( r^2+D+B \varphi^D(r)\Big)^{-\mu/2}  \\
	&\le \frac{1}{c_2}  \Big(B\varphi^D(r)\Big)^{-\mu/2}
	\le  1, \qquad  r\in [0,r_0],
\end{aligned}
\]
because $\varphi^D$ decreases. For large $r$,
we may use (\ref{33.3}) and (\ref{33.6}) to estimate
\[
\begin{aligned}
	v_0(r) ^{-2/\mu}&\le \Big(V_D(r)-cr^{-l}\Big)^{-2/\mu} \\
	&\le V_D(r) ^{-2/\mu}  + \gamma V_D(r) ^{- 2/\mu -1}  r^{-l} \\
	&= r^2+D+ \gamma  \Big( \frac{r^2+D}{r^2} \Big)^{(\mu+2)/2}
r^{-(l-\mu-2)} \\
	&\le r^2+D+\gamma  (D+1)^{(\mu+2)/2}  r^{-(l-\mu-2)}
	\qquad \mbox{for all } r>r_0,
\end{aligned}
\]
whereas 
\[
\begin{aligned}
	\uv(r,0) ^{-2/\mu} &= r^2+D+ B \varphi^D(r) \\
	&\ge  r^2+D + c_1 B r^{-(l-\mu-2)} \qquad \mbox{for all } r>r_0
\end{aligned}
\]
by (\ref{33.4}). Therefore (\ref{33.8}) completes the proof
of (\ref{33.10}), which in turn implies that $v \ge \uv$
holds for all $r\ge 0$ and $t\ge 0$ by comparison.
As a consequence we infer that if we pick $c_4>0$ small such that
$(1+z)^{-\mu/2} \le 1-c_4 z$ for all $z\in [0,B/D]$ then
\bas
	v(r,t) \le (r^2+D)^{-\mu/2}  - c_4 B e^{-\alpha t}  (r^2+D)^{-
(\mu+2)/2}  \varphi^D(r),
	\qquad r\ge 0,\quad t\ge 0,
\eas
which immediately implies (\ref{33.2}) for some $C>0$
 upon invoking (\ref{19.2}).
\qed

\mysection{Optimality}
We proceed to show the optimality of the convergence rates
that we have found in the previous section.

Let us first make sure that convergence {\em from below}
does not occur at a rate faster than asserted by
Lemma~\ref{lem33}.

\begin{lem}\label{lem38}
Suppose that there exist $D>0$, $c>0$ and $l \in (\mu+2,l_\star)$ such that
\be{38.1}
	v_0(r) < V_D(r) \qquad \mbox{for all } r\ge 0
\ee
and
\be{38.2}
	v_0(r) \le V_D(r) - c r^{-l} \qquad \mbox{for all } r\ge 1.
\ee
Then there exists $C>0$ such that the solution $v$ of \eqref{FPeqn}
satisfies
\be{38.3}
	v(r,t) \le V_D(r) - C e^{-(l-\mu-2)(n-l)t}  (r+1)^{-l}
	\qquad \mbox{for all } r\ge 0 \mbox{ and } t\ge 0.
\ee
In particular, $v$ satisfies
\be{38.4}
	\sup_{r\ge 0} \Big(V_D(r)-v(r,t) \Big) \ge C e^{-(l-\mu-2)(n-l)t}
	\qquad \mbox{for all }  t\ge 0.
\ee
\end{lem}

\proof
In view of (\ref{38.1}) and (\ref{38.2}) it is possible to fix a
number $E>D$ close enough to $D$ such that
\be{38.5}
	E-D \le \frac{2\alpha}{\mu(l-\mu-2)^2}
\ee
and
\be{38.6}
	v_0(r) \le V_E(r) \qquad \mbox{for all } r\le 1,
\ee
where again $\alpha:=(l-\mu-2)(n-l)>0$.
Applying Lemma~\ref{lem19}, we obtain a constant $c_1>0$
 such that the solution of (\ref{19.1}) with parameter $d=E+1$
satisfies
\be{38.7}
	\varphi^{E+1}(r) \le c_1 r^{-(l-\mu-2)} \qquad \mbox{for all } r\ge
1,
\ee
and finally we pick $B>0$ so small that
\be{38.8}
	B \le \frac{2c}{\mu c_1 (D+1)^{(\mu+2)/2}}.
\ee
Then
\be{38.9}
	\ov(r,t):=\max \bigg\{ V_E(r), \ \Big(r^2+D+B e^{-\alpha t}
	\varphi^{E+1}(r)\Big)^{-\mu/2}  \bigg\},
	\qquad r\ge 0, \ t\ge 0,
\ee
satisfies
\bas
	\ov(r,0) \ge V_E(r) \ge v_0(r) \qquad \mbox{for all } r\le 1
\eas
by (\ref{38.6}), whereas (\ref{38.2}), (\ref{38.7}), (\ref{38.8})
and the convexity of $0<z\mapsto z^{-2/\mu}$
ensure that
\[
 \begin{aligned}
	\ov(r,0)^{-2/\mu} - v_0(r)^{-2/\mu}
	&\le \Big(r^2+D+B\varphi^{E+1}(r)\Big)
	-\Big( V_D(r)- cr^{-l}\Big)^{-2/\mu} \\
	&\le \Big(r^2+D+B\varphi^{E+1}(r)\Big) \\
	&\qquad -\Big( r^2+D+ \frac{2c}{\mu} (r^2+D)^{-(\mu+2)/2} r^{-l} \Big) \\
	&\le c_1 B r^{-(l-\mu-2)} - \frac{2c}{\mu}
	\Big(\frac{r^2}{r^2+D} \Big)^{(\mu+2)/2}  r^{-(l-\mu-2)} \\
	&\le 0 \qquad \mbox{for all } r\ge 1,
\end{aligned}
\]
because $r^2/(r^2+D)\le 1/(D+1)$ for such $r$. We thus conclude
that $\ov \ge v$ at $t=0$, and next claim that
$\parab \ov \ge 0$ for $r>0$ and $t>0$.
According to Lemma~\ref{lem18} and the fact that $V_E$ is a solution of
(\ref{FPeqn}),
we need to check that for $\psi:=\varphi^{E+1}$,
\[
 \begin{aligned}
	\AD [Be^{-\alpha t}]\psi &= (r^2+D)  \Big( \psi_{rr}+ \frac{n-1}{r}
\psi_r\Big)
	- \mu r \psi_r + \alpha \psi \\
	&  \qquad - B e^{-\alpha t}  \bigg\{ -\psi  \Big( \psi_{rr}+
\frac{n-1}{r} \psi_r \Big)
	+ \frac{\mu}{2} \psi_r^2 \bigg\}
\end{aligned}
\]
is nonnegative whenever $\ov(r,t)>V_E(r)$, that is, when $B e^{-\alpha t}
\psi(r) < E-D$.
At such points, using
(\ref{19.1}) we have
\[
 \begin{aligned}
	\AD [Be^{-\alpha t}]\psi &> (r^2+D)  \Big( \psi_{rr}+ \frac{n-1}{r}
\psi_r\Big)
	- \mu r \psi_r + \alpha \psi \\
	&  \qquad + (E-D)  \Big( \psi_{rr}+ \frac{n-1}{r} \psi_r\Big)
	- \frac{\mu(E-D)\psi_r^2}{2\psi}    \\
	&= (r^2+E)  \Big( \psi_{rr}+ \frac{n-1}{r} \psi_r\Big)
	- \mu r \psi_r + \alpha \psi \\
	& \qquad  - \frac{\mu(E-D)\psi_r^2}{2\psi} \\
	&= -\frac{\mu r}{r^2+E+1}  \psi_r + \frac{\alpha}{r^2+E+1}  \psi
	- \frac{\mu(E-D)\psi_r^2}{2\psi} \\
	&\ge \frac{\alpha}{r^2+E+1}  \psi
	- \frac{\mu(E-D)\psi_r^2}{2\psi},
\end{aligned}
\]
since $\psi$ decreases. Thus, in view of Lemma~\ref{lem19}
and (\ref{38.5}) we infer that
\[
 \begin{aligned}
	\AD [Be^{-\alpha t}]
	\psi &> \bigg\{ \alpha - \frac{\mu(l-\mu-2)^2 (E-D)}{2} \cdot
\frac{r^2}{r^2+E+1} \bigg\}
	 \frac{\psi}{r^2+E+1} \\
	&\ge \bigg\{ \alpha - \frac{\mu(l-\mu-2)^2 (E-D)}{2} \bigg\}
\frac{\psi}{r^2+E+1} \\[1mm]
	&\ge 0
\end{aligned}
\]
is valid whenever $\ov > V_E$, which shows that indeed $\ov$ is
a supersolution of (\ref{FPeqn}), so that the comparison
principle states that $v \le \ov$ in $[0,\infty) \times [0,\infty)$. Since
\[
	r^2+D+B e^{-\alpha t} \psi(r) \le r^2+D+B e^{-\alpha t} \le r^2+E
	\qquad \mbox{for all } r\ge 0
\]
when $t$ is sufficiently large, using (\ref{38.9}) we easily derive
(\ref{38.3}) and (\ref{38.4}).
\qed

\

A similar result on the rate of convergence {\em from above}
can be obtained in a slightly easier manner.

\begin{lem}\label{lem40}
Assume that
\be{40.1}
	v_0(r)>V_D(r) \qquad \mbox{for all } r\ge 0
\ee
and
\be{40.2}
	v_0(r) \ge V_D(r) + cr^{-l} \qquad \mbox{for all } r\ge 1
\ee
are valid with some $D>0$, $c>0$ and $l \in (\mu+2,l_\star)$.
Then there exists $C>0$ such that
\be{40.3}
	v(r,t) \ge V_D(r) + C e^{-(l-\mu-2)(n-l)t}  (r+1)^{-l}
	\qquad \mbox{for all } r\ge 0 \mbox{ and } t\ge 0, 
\ee
and in particular 
\be{40.4}
	\sup_{r\ge 0} \Big( v(r,t)-V_D(r)\Big) \ge C e^{-(l-\mu-2)(n-l)t}
	\qquad \mbox{for all } t\ge 0.
\ee
\end{lem}

\proof
According to (\ref{40.1}) and (\ref{40.2}), it is easy to see that for some
$c_1>0$ we have
\bas
	v_0(r)^{-2/\mu} \le V_D(r)^{-2/\mu}  - c_1 (r^2+D)^{-(l-\mu-2)/2}
\qquad \mbox{for all } r\ge 0.
\eas
In view of Lemma~\ref{lem19}, this means that if $B>0$
 is small enough then taking $\varphi^D$ from (\ref{19.1}) with
$d=D$, we have
\bas
	v_0(r)^{-2/\mu}  \le V_D(r)^{-2/\mu}  - B  \varphi^D(r) \qquad
\mbox{for all } r\ge 0.
\eas
Hence, the function $\uv$ defined by
\bas
	\uv(r,t):=\Big(r^2+D-B e^{-\alpha t}  \varphi^D(r) \Big)^{-\mu/2} ,
\qquad r\ge 0, \ t\ge 0,
\eas
again with $\alpha:=(l-\mu-2)(n-l)>0$, does not exceed $v$ at $t=0$. Since
\[
\begin{aligned}
	\AD[-Be^{-\alpha t}]
	 \varphi^D &= (r^2+D)  \Big( \varphi^D_{ rr} + \frac{n-1}{r}
\varphi^D_{ r} \Big)
	- \mu r \varphi^D_{ r} + \alpha \varphi^D \\
	& \qquad + B e^{-\alpha t }
	 \bigg\{ -\varphi^D  \Big( \varphi^D_{ rr} + \frac{n-1}{r}
\varphi^D_{ r} \Big)
	+ \frac{\mu}{2} (\varphi^D_{ r})^2 \bigg\}
\end{aligned}
\]
satisfies
\[
	\AD [-Be^{-\alpha t}]\varphi^D \ge (r^2+D)
	\Big( \varphi^D_{ rr} + \frac{n-1}{r} \varphi^D_{ r } \Big)
	- \mu r \varphi^D_{ r} + \alpha \varphi^D =0
\]
for all $r>0$ and  $t>0$
by Lemma~\ref{lem19}, we conclude using Lemma~\ref{lem18}
and the comparison principle that $v \ge \uv$ for $r\ge 0$
and $t\ge 0$, which easily results in (\ref{40.3}) and (\ref{40.4}).
\qed

\

\mysection{Universal lower bound for the convergence rate}
{\it Proof of Theorem~\ref{subthm}.} \
 Let us first assume that $0<v_0<V_D$ in $\R^n$ for some $D>0$. Then given $\eps>0$, we write $\alpha:=\alpha_\star+\eps$
 and let $\varphi=\varphi^d(r)$ denote the solution of (\ref{19.1}) corresponding to this value of $\alpha$ and $d:=D+1$.
 Then since $\alpha>\alpha_\star$, Lemma~\ref{le:alphazero} says that there exists $r_0>0$ such that
 $\varphi(r)>0$ for $r\in [0,r_0)$ and $\varphi(r_0)=0$.
 Moreover, from (\ref{19.1}) it is immediately clear that $\varphi_r<0$ in $(0,r_0)$, because at each point in $(0,r_0)$
 where $\varphi_r$ vanishes we must have $\varphi_{rr}<0$ by positivity of $\alpha$ and $\varphi$.
 Since $\varphi_r(r_0)$ cannot vanish in view of an ODE uniqueness argument, it follows that $\varphi_r<0$ even in $(0,r_0]$.
 This in turn entails that
 \be{s11}
	\varphi_{rr} + \frac{n-1}{r} \varphi_r = \frac{\mu r \varphi_r - \alpha \varphi}{r^2+d} \le -c_1
	\qquad \mbox{in } (0,r_0)
 \ee
 is valid for some $c_1>0$, so that by smoothness of $\varphi$,
 \be{s1}
	c_2:=\frac{c_1}{\sup_{r\in (0,r_0)} \Big\{ \varphi \big| \varphi_{rr} +\frac{n-1}{r} \varphi_r\big|
	+\frac{\mu}{2} \varphi_r^2\Big\} }
 \ee
 is positive.

 Now since $v_0(x)<V_D(x)=(|x|^2+D)^{-\frac{\mu}{2}}$ in $\R^n$ and $v_0$ is continuous, it is possible to fix
 $y_0\in (0,c_2]$ small enough such that
 \be{s2}
	v_0(x) \le \Big( |x|^2 + D + y_0 \varphi(|x|) \Big)^{-\frac{\mu}{2}}
	\qquad \mbox{for } |x| \le r_0,
 \ee
 whence the function
 \bas
	\ov(r,t):= \Big(r^2+D+y(t) \varphi(r)\Big)^{-\frac{\mu}{2}}, \qquad (r,t) \in [0,r_0] \times [0,\infty),
 \eas
 with
 \bas
	y(t):=y_0 e^{-\alpha t}, \qquad t\ge 0,
 \eas
 satisfies
 \be{s3}
	v_0(x) \le \ov(|x|,0) \qquad \mbox{for } |x| \le r_0.
 \ee
 Since clearly $v_0<V_D$ in $\R^n$ implies $v(x,t) \le V_D(x)$ in $\R^n \times (0,\infty)$ by comparison, we also have
 \be{s4}
	v(x,t) \le V_D(x) = \ov(|x|,t)
	\qquad \mbox{whenever $|x|=r_0$ and } t\ge 0.
 \ee
 Furthermore, recalling the definition (\ref{2.100}) of $\parab$, by Lemma~\ref{lem18} we have
 \be{s5}
	\parab \ov = \frac{\mu}{2} y(t) \Big( r^2+D+y(t) \varphi(r)\Big)^{-\frac{\mu+2}{2}} \AD[y(t)] \varphi
	\qquad \mbox{for $r\in(0,r_0)$ and } t>0,
 \ee
 where
 \bas
	\AD[y(t)] \varphi &=& (r^2+D)\Big(\varphi_{rr}+\frac{n-1}{r}\varphi_r\Big) - \mu r \varphi_r - \frac{y'}{y} \varphi \\
	& & \hspace*{10mm}
	- y \Big\{ -\varphi \Big(\varphi_{rr}+\frac{n-1}{r}\varphi_r\Big) + \frac{\mu}{2}\varphi_r^2 \Big\} \\
	&=& (D-d) \Big(\varphi_{rr}+\frac{n-1}{r}\varphi_r\Big) + \Big( - \frac{y'}{y} - \alpha \Big) \varphi \\
	& & \hspace*{10mm}
	- y \Big\{ -\varphi \Big(\varphi_{rr}+\frac{n-1}{r}\varphi_r\Big) + \frac{\mu}{2}\varphi_r^2 \Big\}
	\quad \mbox{for $r\in(0,r_0)$ and } t>0
 \eas
 according to (\ref{19.1}). Since $d=D+1$ and $y'/y \equiv -\alpha$, we thus have
 \bas
	\AD[y(t)] \varphi &=& - \Big(\varphi_{rr} +\frac{n-1}{r}\varphi_r\Big)
	- y(t) \Big\{ -\varphi \Big(\varphi_{rr}+\frac{n-1}{r}\varphi_r\Big) + \frac{\mu}{2}\varphi_r^2 \Big\} \\
	&\ge& c_1 - y_0 \Big\{ \varphi \Big|\varphi_{rr}+\frac{n-1}{r}\varphi_r\Big| + \frac{\mu}{2}\varphi_r^2 \Big\} \\
	&\ge& 0
	\qquad \mbox{for $r\in(0,r_0)$ and } t>0
 \eas
 because of (\ref{s11}), (\ref{s1}) and the fact that $y_0 \le c_2$.
 In view of (\ref{s3})-(\ref{s5})  and  a comparison argument this implies that
 \bas
	v(x,t) \le \ov(|x|,t) \qquad \mbox{for $|x| \le r_0$ and } t\ge 0,
 \eas
 which in particular means that
 \bas
	V_D(0)-v(0,t) \ge D^{-\frac{\mu}{2}} - (D+y_0 e^{-\alpha t} )^{-\frac{\mu}{2}}
	\qquad \mbox{for all } t\ge 0.
 \eas
 Now an application of the elementary inequality
 \[(1+z)^{-\frac{\mu}{2}} \le 1-\frac{\mu}{2}(1+z_0)^{-\frac{\mu+2}{2}} z,
\qquad
0\le z \le z_0 < \infty,
\]
 to $z:=z_0 e^{-\alpha t}$ and $z_0:=y_0/D$ shows that
 \bas
	V_D(0)-v(0,t) \ge D^{-\frac{\mu}{2}} \frac{\mu}{2} \Big(1+\frac{y_0}{D}\Big)^{-\frac{\mu+2}{2}} e^{-\alpha t}
	\qquad \mbox{for all } t\ge 0
 \eas
 and thereby establishes (\ref{1.2.1}) in the case when $0<v_0<V_D$.

 The proof in the case $v_0>V_D$ can be carried out along the same lines.   Indeed,
 \bas
	\uv(r,t):=\Big(r^2+D-y(t) \varphi(r)\Big)^{-\frac{\mu}{2}}, \qquad (r,t) \in [0,r_0] \times [0,\infty),
 \eas
 becomes a subsolution of (\ref{FPeqn}) satisfying $v(x,t) \ge \uv(|x|,t)$ in $\R^n \times (0,\infty)$, provided that
 $r_0$ and $\varphi$ are chosen as above and $y(t):=\tilde y_0 e^{-\alpha t}$, where $\alpha$ is as before and $\tilde y_0>0$
 is small enough fulfilling
 \bas
	v_0(x) \ge \Big(|x|^2 + D - \tilde y_0 \varphi(|x|) \Big)^{-\frac{\mu}{2}} \qquad \mbox{for } |x| \le r_0
 \eas
 and $\tilde y_0 < \min\{c_2,D\}$ with $c_2$ as given by (\ref{s1}).
 In conclusion, in this case we obtain
 \bas
	v(0,t)-V_D(0) &\ge& \Big(D-\tilde y_0 e^{-\alpha t}\Big)^{-\frac{\mu}{2}} - D^{-\frac{\mu}{2}} \\
	&\ge& D^{-\frac{\mu}{2}} \frac{\mu}{2} e^{-\alpha t}
	\qquad \mbox{for all } t\ge 0
 \eas
 by convexity of $[0,1) \ni z \mapsto (1+z)^{-\frac{\mu}{2}}$. This again proves (\ref{1.2.1}).
\qed

\section{Comments}

 We have not discussed the problem of constructing solutions in the 
non-variational basin for $m\ge m_*$. For $m> m_*$, the existence of such a
basin seems unlikely but for $m= m_*$ we expect that the variational basin
studied in \cite{BGV} represents just a part of the whole basin of
attraction and the variational rate found in \cite{BGV} is the fastest
possible rate of convergence. The construction of solutions with slower rates
for $m= m_*$ seems challenging and it will become a theme of future
investigations.

 The questions of extinction behaviour are different when the problem 
is posed in a bounded domain. Then the extinction rate for bounded solutions  
is universal, of the form
$\|u(\cdot,\tau)\|_\infty = O((T-\tau)^{1/(1-m)})$ when $m>
m_s:=(n-2)/(n+2)$, cf. \cite{BH},  but the question is more complicated when
$m\le m_s$. Convergence rates for $m$ near 1 have been recently obtained in \cite{BGV10}.

\

\noindent
{\bf Acknowledgments.}
M.~Fila was supported by the Slovak Research and Development Agency
under the contract No. APVV-0134-10 and by VEGA grant 1/0465/09 (Slovakia).
E.~Yanagida was supported  by the Grant-in-Aid for Scientific
Research (A)  (No.~19204014) from the
Ministry of Education, Science, Sports and Culture (Japan).
J.~L.~V\'azquez was supported by Project MTM2008-06326 (Spain).
 
\newpage

%

\


\noindent{\sc Addresses}
\bigskip

\noindent Marek Fila\\
Department of Applied Mathematics and Statistics, Comenius University,\\
84248 Bratislava, Slovakia

\medskip
\noindent Juan Luis V\'azquez\\
Departamento de Matem\'aticas, Universidad Aut\'onoma de Madrid,\\
28049 Madrid, Spain

\medskip
\noindent Michael Winkler \\
Institut f\"ur Mathematik, Universit\"at Paderborn,\\
33098 Paderborn, Germany

\medskip
\noindent Eiji Yanagida \\
Department of Mathematics, Tokyo Institute of Technology,\\
Meguro-ku, Tokyo 152-8551, Japan
\end{document}